\documentclass[12pt]{article}
\usepackage{multirow,verbatim}
\usepackage{latexsym, amsbsy, amssymb, geometry}
\usepackage{natbib}
\usepackage{lscape}
\usepackage{float}
\usepackage{booktabs}
\usepackage{graphicx}
\usepackage{algorithm}
\usepackage{algpseudocode}
\renewcommand{\baselinestretch}{1.6}

\usepackage{tikz}
\usepackage{graphicx}
\usepackage{epsfig,float}
\usepackage{amsmath}
\usepackage{mathrsfs}
\usepackage{array}
\usepackage{enumerate} 
\usepackage[shortlabels]{enumitem}
\usepackage{multirow}
\usepackage{latexsym, amsbsy, amssymb}
\usepackage{verbatim}
\usepackage{color}
\usepackage{url}
\usepackage{caption}
\usepackage{subcaption}

\usepackage{xr}
\usepackage{xcolor, color}

\newtheorem{assumption}{{\bf Assumption}}
\newtheorem{remark}{Remark}
\newtheorem{theorem}{Theorem}
\newtheorem{lemma}[theorem]{Lemma}

\newcommand{\indep}{\;\, \rule[0em]{.03em}{.6em} \hspace{-.25em}
	\rule[0em]{.65em}{.03em} \hspace{-.25em}
	\rule[0em]{.03em}{.6em}\;\,}

\newcommand{\trans}{^{\mbox{\tiny {\sf T}}}}

\def\nano{\scriptscriptstyle}
\newcommand\hi[1]{^{\nano #1}}
\def\real{\mathbb R}
\newcommand\ca[1]{{\cal{#1}}}
\newcommand\lo[1]{_{\nano #1}}

\def\nano{\scriptscriptstyle}

\def\nano{\scriptscriptstyle}

\def\ka{\kappa}

\newcommand{\RN}[1]{
	\textup{\uppercase\expandafter{\romannumeral#1}}
}

\newcommand{\ben}{\begin{enumerate}}
	\newcommand{\een}{\end{enumerate}}
\def \mbb {\mathbb}

\def\shortbar{\rule[-.0em]{.04em}{.15em}}
\def\bbar{\mbox{\raisebox{-.05em}{$\,\overset{\overset{\shortbar}{\shortbar}}{\shortbar}\,$}}}

\begin{document}

	\begin{center}
		{\Large{\bf
			Variable Selection for Additive Global Fr\'echet Regression}}

        \vskip.5cm
			{\large 
            Haoyi Yang$^\ast$, Satarupa Bhattacharjee$^\dagger$, Lingzhou Xue$^\ast$, and Bing Li$^\ast$
            }\\
		\vskip.3cm
			$^\ast$Department of Statistics, The Pennsylvania State University\\
            $^\dagger$Department of Statistics, University of Florida\\

        \vskip.5cm
        First Version: October 2024\\
        This Version: August 2025

   	\end{center}
	
	\begin{abstract}
We present a novel framework for variable selection in Fr\'echet regression with responses in general metric spaces, a setting increasingly relevant for analyzing non-Euclidean data such as probability distributions and covariance matrices. Building on the concept of (weak) Fr\'echet conditional means, we develop an additive regression model that represents the metric-based discrepancy of the response as a sum of covariate-specific nonlinear functions in reproducing kernel Hilbert spaces (RKHS). To address the absence of linear structure in the response space, we transform the response via squared distances, enabling an interpretable and tractable additive decomposition. Variable selection is performed using Elastic Net regularization, extended to the RKHS setting, and further refined through a local linear approximation scheme that incorporates folded concave penalties such as the SCAD. We establish theoretical guarantees, including variable selection consistency and the strong oracle property, under minimal assumptions tailored to metric-space-valued responses. Simulations and applications to distributional and matrix-valued data demonstrate the scalability, interpretability, and practical effectiveness of the proposed approach. This work provides a principled foundation for statistical learning with random object data.

    \end{abstract}


\section{Introduction}
\label{sec:intro}

Recent years have witnessed a rapid growth of complex data types that lie outside traditional Euclidean spaces, such as probability distributions \citep{panaretos2019statistical, petersen2022modeling, chen2023wasserstein, lin2023causal, zhang2024nonlinear,bhattacharjee2025doubly}, covariance matrices \citep{dryden2009non, lin2019riemannian}, 
and trajectories or curves \citep{bhattacharjee2025geodesic, dubey2020functional}. Often referred to as \emph{random objects}, such data reside in metric spaces, where arithmetic operations such as addition or scalar multiplication are not defined. As a result, classical regression methods are inadequate, motivating statistical frameworks that explicitly account for the intrinsic geometry of these spaces.

Within this context,  Fr\'echet regression has emerged as a key approach for modeling metric-space-valued responses with Euclidean or RKHS-valued covariates~\citep{pete:19, bhat:25}. Recent work has advanced estimation and inference for these models \citep{lin2019intrinsic, bhattacharjee2025geodesic, kim2020nonparametric}, but the issue of variable selection remains underexplored. Identifying a sparse and informative subset of covariates is crucial for enhancing interpretability and alleviating the curse of dimensionality, which poses a major challenge in high-dimensional and functional object regression \citep{bhattacharjee2023single, zhang2024dimension}. A principled selection strategy can improve both predictive performance and scientific insight.

In Euclidean regression, variable selection is supported by a rich literature. The Lasso~\citep{tibshirani1996regression} introduced an $\ell_1$-penalty to induce sparsity, while the Elastic Net~\citep{zou2005regularization} combined $\ell_1$ and $\ell_2$ penalties  to improve stability under collinearity. Nonconvex penalties such as the SCAD~\citep{fan2001variable} and adaptive weighting schemes such as Adaptive Lasso~\citep{zou2006adaptive} further improve selection consistency and can achieve the oracle property. 

Additive models provide a flexible framework for nonlinear regression and variable selection, where the conditional expectation is expressed as
$
\mbb{E}[Y|X_1\dots,X_p] = \sum_{j=1}^p f_j(X_j),
$
with each $f_j$ capturing the marginal effect of $X_j$. Originating from \citet{stone1985additive} and extended to functional data by \citet{muller2008functional}, this structure enables univariate smoothing and simplifies multivariate estimation. Penalization approaches such as the group Lasso and component-wise penalties~\citep{meier2009high, ravikumar2009sparse} adapt the Lasso to the additive setting, providing variable selection consistency while preserving interpretability.

However, these Euclidean approaches for variable selection and additive modeling do not naturally extend to general metric-space-valued responses due to the absence of algebraic structure (such as addition and multiplication). \citet{tucker2019genomics} studied variable selection in a functional setting, but their framework is not directly applicable to general Fréchet regression. More recently, \citet{tucker2023variable} proposed the Fr\'echet Ridge Selection Operator (FRiSO) for global linear Fr\'echet regression, extending individually penalized ridge regression \citep{wu2021can} from the Euclidean case. Their method assigns a separate ridge penalty to each predictor, optimizes these penalties under an $\ell_1$-type budget constraint, and then identifies active variables based on whether the corresponding ridge penalty is zero. In practice, as implemented in their code, variable selection is carried out by post-thresholding the ridge penalties, which is highly sensitive to the chosen threshold(shown in Section~\ref{suppl:sec:friso_refit} of the supplements), as we illustrate in the supplement. On the theoretical side, Theorem 1 of \citet{tucker2023variable} shows only that the penalties for active variables converge to zero and those for inactive variables diverge to infinity in probability, but, as noted in their remark after Theorem 1, this does not establish the desired almost-sure selection consistency.

Motivated by these gaps, we propose a new framework that transplants additive modeling from objects to distances, thereby enabling a unified approach for both linear and nonlinear Fr\'echet regression, avoiding the need for ad-hoc post-thresholding, and delivering rigorous variable selection guarantees. Specifically, we transform the response via squared distances as $u(y) = d_Y^2(Y,y) - d^2(Y,y_0)$, where $y_0$ is a fixed reference point (e.g., the Fr\'echet mean), mapping the response into a scalar function amenable to additive modeling.

Our framework models this transformed response as an additive function of the covariates, with the nonlinear component of each covariate represented flexibly using reproducing kernel Hilbert spaces (RKHS). This structure decomposes the metric relationship between the response and covariates into interpretable, covariate-specific effects. To perform variable selection, we formulate a penalized RKHS optimization problem, extending Elastic Net regularization~\citep{zou2005regularization} to the non-Euclidean and nonlinear setting. This approach promotes both sparsity and stability while respecting the geometry of the response space, resulting in a parsimonious and interpretable model that applies to complex data types.

Theoretically, we establish a fundamental \emph{sparsity invariance} property in Theorem \ref{theorem: Negative} for a broad class of additive models defined in the negative-type metric space, including many linear model examples in \citet{tucker2023variable}. This property ensures that the active set of predictors remains unchanged regardless of the choice of $y_0$ in $u(y)= d_Y^2(Y,y) - d^2(Y,y_0)$. Motivated by this insight, we formalize Assumption \ref{ass:sparsity:invariance} as the key condition underpinning the interpretability and stability of variable selection in metric spaces. Building on the sparsity invariance, and without imposing additional assumptions, our analysis establishes variable selection consistency and the strong oracle property for both linear and nonlinear Fr\'echet regression models, which are not previously available in the literature. To derive these guarantees, we also extend several core analytical tools from Euclidean settings to reproducing kernel Hilbert spaces (RKHS): in particular, we generalize the Karush–Kuhn–Tucker (KKT) conditions and develop probabilistic inequalities that remain valid in infinite-dimensional RKHS settings. These technical contributions are central to our proofs and may also be of independent interest for future research in metric-space regression.

The remainder of the paper is structured as follows: Section 2 reviews background on Fr\'echet regression and RKHS tools. Section 3 introduces our modeling framework, response transformation, and additive structure. Section 4 presents the regularization scheme and implementation. Section 5 provides theoretical results on consistency and sparsity recovery. Section 6 illustrates the method via simulations and data applications. Section 7 concludes. The complete proofs and additional numerical results are presented in the supplement.

\section{Preliminaries} \label{sec:prelim}

\textbf{Notations and Setup.} Let $(\Omega, \ca F, P)$ be a probability space. Let $(\Omega_X, d_X)$ and $(\Omega_Y, d_Y)$ be metric spaces, where $d_X$ and $d_Y$ are the metrics defined on the sets $\Omega_X$ and $\Omega_Y$, respectively.  Let $\ca F_X$ and $\ca F_Y$ be the Borel $\sigma$-fields in $\Omega_X$ and $\Omega_Y$ corresponding to the open sets determined by $d_X$ and $d_Y$.  Let $X: \Omega \to \Omega_X$ and $Y: \Omega \to \Omega_Y$ be random elements that are measurable, respectively, with respect to $\ca F/ \ca F_X$ and $\ca F / \ca F_Y$.  Such random elements are called {\em statistical objects}. 
In our framework, $\Omega\lo X$ is a product metric space, that is, $X = (X_1, \ldots, X_p)$ is a $p$-dimensional random object covariates, with each $X_j : \Omega \to \Omega_{X_j}$, where $(\Omega_{X_j}, d_{X_j})$ is a separable metric space for each $j=1,\dots,p$, where $p\geq 1$ can be moderately large but not increasing with the sample size. Throughout, we assume that the random elements $(X, Y)$ are jointly measurable and that $P$ is the joint distribution of $(X, Y)$ on $\prod_{j=1}^p \Omega_{X_j} \times \Omega_Y$.

We assume the existence of a positive definite kernel $\kappa_j: \Omega_{X_j} \times \Omega_{X_j} \to \real$ for the $j$-th covariate, and let $\ca M\lo j $ be the RKHS generated by $\ka_j$. That is, for any finite set $\{x_1, \ldots, x_m\} \subset \Omega_{X_j}$, the matrix $[\kappa_j(x_i, x_k)]_{1 \leq i,k \leq m}$ is symmetric and positive semi-definite. The RKHS $\mathcal{M}_j$ is the unique Hilbert space of functions $f_j: \Omega_{X_j} \to \mathbb{R}$ such that: (i) For any $x \in \Omega_{X_j}$, the function $\kappa_j(\cdot, x)$ belongs to $\mathcal{M}_j$; (ii) For any $f_j \in \mathcal{M}_j$ and $x \in \Omega_{X_j}$, the reproducing property holds: $f_j(x) = \langle f_j, \kappa_j(\cdot, x) \rangle_{\mathcal{M}_j}$. While there are sufficient conditions for a metric space to admit such kernels, we simply impose this as a standing assumption.

\begin{assumption}\label{ass:positive kernel} 
There is a positive definite kernel $\ka_X: \Omega_X \times \Omega_X \to \real$.
\end{assumption}
If $\Omega_X$ is of negative type, the metric-induced kernel is positive definite \citep{sejdinovic2012hypothesis}. \citet{zhang2024dimension} showed that Gaussian and Laplacian kernels defines as $\kappa_G(x,x') = \exp(-\gamma\lo X d\lo X \hi 2(x, x'))$ and {$\kappa_L(x, x') = \exp(-\gamma_X d\lo X(x, x'))$}, respectively are positive definite on complete, separable metric spaces, and their RKHS is dense in $L^2(P_X)$. These kernels perform well for common covariates such as univariate distributions with Wasserstein-2 distance and symmetric positive definite matrices with Frobenius distance.  

Embedding each covariate into an RKHS enables flexible modeling of object-valued regression by accommodating a spectrum of model complexities, such as the linear model, the polynomial model, and a family of functions that is dense in the $\ell_2$ space. In particular, our primary focus is variable selection in object regression, for which the additive modeling framework via RKHS embedding of the covariate metric spaces facilitates an interpretable framework. 
No embedding is assumed for $\Omega_Y$. Indeed, any function that fits $Y$ serves our purpose, and embedding is just one way to generate such a function for additive models. Our embedding process does not disrupt the interpretation of the random object, as the interpretation remains consistent regardless of the function applied to it.

\textbf{Fr\'echet Regression.} For random object responses $Y$, expectation and variance have been generalized to metric spaces via the Fr\'echet mean and variance~\citep{frec:48}: $\mathbb{E}\lo \oplus(Y) = \arg\min_{y \in \Omega_Y} \mathbb{E}[d_Y^2(Y, y)]$, and $V_\oplus = \mathbb{E}[d_Y^2(Y, \mathbb{E}\lo \oplus(Y))]$. To model regression between $Y \in (\Omega\lo Y, d\lo Y)$ and covariate $X \in \real^p$, \citet{pete:19} proposed the Fr\'echet regression:
\begin{align}
    \mathbb{E}\lo \oplus (Y|X) &= \arg\min_{y \in \Omega_Y} \mathbb{E}[d_Y^2(Y, y) \mid X] 
\end{align}
which interprets conditional expectation through conditional Fr\'echet means. Estimation under this framework entails modeling the joint distribution of $(Y, X)$, by viewing the regression function as an alternative target of weighted Fr\'echet means, with weights that change globally linearly (or locally) with the covariates and are derived from those of the corresponding standard multiple linear regression (or local linear kernel regression) with Euclidean responses. While these methods provide useful analogs, the global model can be overly restrictive, and local methods suffer from the curse of dimensionality.

To address these limitations, \citet{bhat:25} proposed a nonlinear global Fr\'echet regression using weak conditional means. 
Intuitively, the weak conditional expectation in a Hilbert space is a generalization of conditional expectation in the sense that, while the latter is the $\ell_2(P_X)$ space projection of the response, the former is a projection to an arbitrary subset (of $\ell_2(P_X)$) of functions situated in an RKHS. In the same spirit, the nonlinear global Fr\'echet regression, $\mbb{E} \lo\oplus(Y\bbar X)$ is the best approximation of $d\hi 2(Y,y)$ among all of functions in $\overline{\ca M}\lo X \subseteq L_2(P\lo X)$, where $\overline{\ca M}\lo X$ is the $\ell_2 (P_X)$-closure of $\ca M_X$. Formally, the weak conditional Fr\'echet mean is the minimizer of the inducing function of the Carleman regression operator $R_{XU(y)}$, which captures the dependence of $d\hi 2(Y,y)$ on $X$.

The nonlinear global Fr\'echet regression via weak Fr\'echet conditional means offers a broad spectrum of flexible models for random object regression. Among many appealing intuitions that it provides as a Fr\'echet regression approach, we pick our motivation from the following aspect:  the scalar function $ E [d \hi 2 (Y, y)| X]$, analogous to classical models:
$ E [d \hi 2 (Y, y)| X]$,  {\em can be} approximated by a basis expansion; that is, it is reasonable to expect
\begin{align*}
E [d \hi 2 (Y, y)| X] \approx c_1 f_1 (X) + \cdots + c_m f_m (X), 
\end{align*}
where $f_j$ are basis functions and $c_j$ are coefficients. While such an expansion is infeasible for the object-valued quantity $E_\oplus(Y|X)$ due to the lack of linear operations in general metric spaces, it is well-defined for the expected squared distance, which is scalar-valued. Thus, the weak Fr\'echet mean selects the $y\in \Omega\lo Y$ minimizing the best approximation to this conditional expectation over a given function class. 
This approach provides a principled way to extend basis expansions and additive models to the non-Euclidean setting of random object regression by modeling the metric-based distances rather than directly modeling the relationship between covariates and response through coefficients.

\section{Methodology}\label{sec:add:fr:reg}

This section first introduces an additive modeling framework in Subsection 3.1 and sparsity invariance in Subsection 3.2, and then present the additive global Fr\'echet regression in Subsection 3.3 and penalized estimation via local linear approximation in Subsection 3.4.

\subsection{Additive Modeling}\label{subsec:additive:fr:model}

Additive regression modeling offers a flexible, interpretable framework that can capture nonlinear relationships between covariates and the response while mitigating the risk of overfitting associated with more complex, fully nonparametric methods. We start by considering a naive approach to the additive modeling of the random object responses with moderately high-dimensional covariates by phrasing the conditional expectation   $\mathbb{E}[d_Y^2(Y, y) \mid X = x] $ as
\begin{equation}
    \mathbb{E}[d_Y^2(Y, y) \mid X = x] \approx \sum_{j=1}^p f_j(X_j),
\end{equation}
for the component functions $f_j,\ j =1,\dots, p$, residing in some suitably rich function class.
However, this equation does not hold even for the linear regression given by $Y=X^T\beta+\epsilon$. We can write 
$
\mathbb{E}[(Y-y)^2\mid X = x]=(X^T\beta-y)^2 = y^2 + (X^T\beta)^2 + 2yX^T\beta,
$
which contains the interaction term, such as $X_jX_k$. However, if we modify the objective function to be 
\[
\mathbb{E}[(Y-y)^2-(Y - y \lo 0)^2\mid X = x]=2(y \lo 0 - y \lo 1) \beta \trans X + y \hi 2 - y \lo 0 \hi 2, 
\]
for any \( y_0 \in \Omega_Y \), then the right-hand side is additive in $X \lo 1, \ldots, X \lo p$. 
This motivates us to model $ \mathbb{E}[d_Y^2(Y, y)- \mathbb{E}[d_Y^2(Y, y_0) \mid X = x],$  for a fixed $y \lo 0 \in\Omega \lo X$,   via an additive structure, which is more reasonable and interpretable.
We define the conditional expected discrepancy as:
\[
D(y, x) := \mathbb{E}[(d_Y^2(Y, y) - d_Y^2(Y, y_0)) \mid X = x],
\]
where \( y_0 \in \Omega_Y \) is a fixed reference point. We define the global Fr\'echet additive model as
\begin{align}
    \label{gl:fr:add: model}
    d_\oplus(x) & = \arg\min_{y \in \Omega_Y} \mathbb{E}[d_Y^2(Y, y) - d_Y^2(Y, y_0) \mid X = x]
\end{align}
Note that, for a conventional additive regression modeling of Hilbert space valued response and covariates pairs $(X,Y)$ of the form $Y = \sum_{j=1}^p f_j(X_j) +\text{error}$, the conditional expected discrepancy $D(y,x)$, which boils down to a difference in Hilbert space norms, can be decomposed as $D(y,x)= \mathbb{E}[\|Y -y\|^2\lo{\ca H} - \|Y-y_0\|^2\lo{\ca H} ] = c_0(y) + \sum_{j=1}^p c_{yj}(x_j)$, where $c_0(y)$ is independent of $x$ and each $c_{yj}(x_j)$ represents the contribution of the covariate $X_j$.

We aim to develop robust variable selection methods for the additive object regression model. Posed as a minimization problem \eqref{gl:fr:add: model}, it is computationally intensive to rerun the variable selection algorithms repeatedly while searching across the space $\{y\in \Omega\lo Y\}.$ However, a critical observation regarding the stability of the additive modeling can be utilized. 

\subsection{Sparsity Invariance}\label{subsec:fund:ass}

For a negative-type metric space, the metric $d\lo Y$ admits an isometric embedding into a Hilbert space. Under such embedding, an additive modeling of the conditional expected discrepancy $D(y,x)$ enjoys the sparsity invariance property, which ensures that once a component function \( f_j \) is inactive for some \( y \), it remains inactive for most other \( y \), allowing consistent variable selection. Thus, the sparsity invariance property can be extended to such metric spaces; the following theorem demonstrates the sparsity structure of the additive model.

\begin{theorem}[Sparsity Invariance in Negative Type Spaces]\label{theorem: Negative}
Suppose that \( (\Omega_Y, d_Y) \) is a negative-type metric space and that after embedding $Y$ into a Hilbert space with mapping $\phi$, we have an additive decomposition for $\phi(Y)$ as $\phi(Y)=\sum_{j=1}^p f_j(X_j)+\epsilon$, where $f_j$ is a mapping from covariate space to the transformed Hilbert space and $\epsilon$ is independent of covariates. Then the conditional expected discrepancy also admits the additive
decomposition:
$\mathbb{E}[d_Y^2(Y, y) - d_Y^2(Y, y_0) \mid X = x] = c_0(y) + \sum_{j=1}^p c_{yj}(x_j)$. Moreover, \( c_{yj}(\cdot)\) is a zero function whenever  
\( \langle \phi(y)-\phi(y_0), f_j(X_j) \rangle = 0 \), and the sparsity pattern is invariant for all \( y \), except on a union of \( p \) subspaces in the embedding space.
\end{theorem}

While sufficient conditions exist for a metric space to possess such a sparsity invariance property, we make this requirement as Assumption \ref{ass:sparsity:invariance}, especially to facilitate computational traceability of the proposed variable selection methods in the subsequent sections. 

\begin{assumption} (Fundamental assumption for stability in variable selection)
\label{ass:sparsity:invariance}
   Let $A_y$ be the smallest subset of $\{1,\dots,p\}$ such that $d(Y,y) \indep X|X_{A_y}$. Then, $A_y = A_{y'}$ for all $y,y'\in \Omega_Y.$
\end{assumption}

This assumption ensures that the active set of covariates remains invariant across all $y\in \Omega\lo Y$, thus fully characterizing variable relevance without requiring an additive model. This assumption is crucial: without 
it, variable selection in general metric spaces would lose 
interpretability, since the active set could vary arbitrarily with the choice of $y$. Nevertheless, additive modeling offers a practical means to disentangle covariate effects, making it well-suited for applying variable selection methods. Under this assumption, the sparsity invariance property holds in strong negative type metric spaces, as established in Theorem~\ref{theorem: Negative}.

Unlike Euclidean settings, the nonlinearity of metric-space-valued responses necessitates modeling through conditional Fr\'echet means. The sparsity invariance assumption thus plays a central role in our framework, enabling the recovery of a stable, globally relevant set of covariates without invoking global model parameters. When this property fails, sparsity may become location-dependent, varying across regions of the metric space. This undermines the notion of a single interpretable active set and complicates model interpretation, highlighting the fundamental importance of sparsity invariance in object regression.

\subsection{Additive Global Fr\'echet Regression}
\label{subsec:gl:add:fr:reg}
In this section, we develop a regularization framework for simultaneous estimation and variable selection. Having expressed the regression relationship between the random object pairs $(X,Y)$ only in terms of the metric $d\hi2(Y,y)$, we can rephrase~\eqref{gl:fr:add: model} as 
\begin{align}
    \label{gl:fr:add: model2}
    d_\oplus(x) & = \arg\min_{y \in \Omega_Y} \mathbb{E}[U(y)], \text{ where }
    U(y) = d_Y^2(Y, y) - d_Y^2(Y, y_0) 
\end{align}
Further, under the sparsity invariance assumption, the minimization problem above is invariant with respect to $y$, thus reducing the computational cost and complexity by a significant amount. Under this assumption, we can reduce model~\eqref{gl:fr:add: model2} as $$d_Y^2(Y, y) - d_Y^2(Y, y_0) =  \sum_{j=1}^p f_j(X_j) + \epsilon,$$ where each component function $f_j$ belongs to the RKHS $\mathcal{M}_j$ associated with the $j$-th covariate, and $\epsilon$ is a random error term satisfying $\mathbb{E}[\epsilon \mid X] = 0$.

The additive structure simplifies modeling by reducing interaction complexity. Joint kernels can capture interactions, but searching over all subsets is combinatorially expensive in high dimensions. Thus, joint modeling is suitable when interactions are known or in low dimensions, but the additive model remains a practical and interpretable choice.

We assume that each $X_j$ lies in a structured function space $\ca M_j$, allowing nonlinear effects of $U(y)$ on $X_j$ to be modeled separately. By the reproducing property, each $f_j\in \ca M_j$ satisfies $f_j(x) = \langle f_j, \kappa_j(\cdot, x) \rangle_{\mathcal{M}_j},$ for $x \in \Omega_{X_j}.$ Evaluated at sample points $\{X_{1j}, \ldots, X_{nj}\}$, $f_j$  admits the kernel expansion: $$f_j = \sum_{k=1}^{\infty} \alpha_{jk} \kappa_j(\cdot, X_{jk}),$$ with coefficients $\{\alpha_{jk}\}$ depending on $f_j$ and $X_{jk}$.
In practice, estimation uses a finite-dimensional approximation based on the data. The RKHS norm is computed as $\|f_j\|_{\mathcal{M}_j}^2 = \sum_{k=1}^{\infty} \sum_{l=1}^{\infty} \alpha_{jk} \alpha_{jl} \kappa_j(X_{jk}, X_{jl})$, providing a tractable regularization framework 
via the kernel.

\subsection{Adaptive Penalization and Local Linear Approximation}
\label{subsec:add:pen:lla}

We estimate the component functions $\{f_j\}_{j=1}^p$ by solving an Elastic Net regularized least squares that promotes sparsity and smoothness. The optimization problem is formulated as
\begin{equation}\label{optprob}
    \min_{f_1, \ldots, f_p} \frac{1}{2n} \sum_{i=1}^n \left( U_i(y) - \sum_{j=1}^p \langle f_j, \kappa_j(\cdot, X_{ij}) \rangle_{\mathcal{M}_j} \right)^2 + \lambda_1 \sum_{j=1}^p \|f_j\|_{\mathcal{M}_j} + \frac{\lambda_2}{2} \sum_{j=1}^p \|f_j\|_{\mathcal{M}_j}^2,
\end{equation}
where
$U_i(y) = d_Y^2(Y_i, y) - d_Y^2(Y_i, y_0)$ is the transformed response for the $i$-th observation; $\langle f_j, \kappa_j(\cdot, X_{ij}) \rangle_{\mathcal{M}_j}$ represents the evaluation of $f_j$ at $X_{ij}$ via the reproducing property; $\|f_j\|_{\mathcal{M}_j}$ is the RKHS norm of $f_j$, serving as a smoothness measure; $\lambda_1 > 0$ controls the sparsity of the component functions by encouraging $\|f_j\|_{\mathcal{M}_j} = 0$ for irrelevant covariates; $\lambda_2 > 0$ penalizes the complexity of nonzero functions to prevent overfitting.

The combination of $\ell_1$ and $\ell_2$ penalties enables simultaneous variable selection and function estimation: the $\ell_1$ term induces sparsity by shrinking some component functions to zero, while the $\ell_2$ term stabilizes estimation by controlling smoothness. Beyond convex penalties, the folded concave penalties like SCAD~\citep{fan2001variable} and MCP~\citep{zhang2010nearly} provide stronger theoretical properties, achieving the strong oracle property \citep{fan2014strong}. Let \( P_{\lambda_1}(\cdot) \) be a general folded concave penalty satisfying the follwoing conditions: 
\begin{assumption}
\label{ass:concave:penalty}
\item[(i)] \( P_\lambda(t) \) is increasing and concave on \( [0, \infty) \), with \( P_\lambda(0) = 0 \);
    \item[(ii)] \( P'_\lambda(t) \) exists for \( t > 0 \), and \( P'_\lambda(0^+) \ge a_1 \lambda \) for some constant \( a_1 > 0 \);
    \item[(iii)] \( P'_\lambda(t) \ge a_1 \lambda \) for all \( t \in (0, a_2 \lambda] \), with \( a_2 > 0 \);
    \item[(iv)] \( P'_\lambda(t) = 0 \) for all \( t \ge a \lambda \), where \( a > a_2 \).
\end{assumption}
\noindent Notable examples of such penalties include SCAD and MCP. The SCAD has derivative: $P'_{\lambda,a}(t) =\lambda$ if $t \leq \lambda,$ $=\frac{a\lambda - t}{a - 1}$ if $\lambda < t \leq a\lambda$, and $=0$ if $t > a\lambda$, 
where $t > 0$ and $a > 2$.

A related variant, SCAD-L$_2$, was introduced by \cite{zeng2014group} as an extension of the Elastic Net, replacing the $\ell_1$ penalty with SCAD while retaining the $\ell_2$ penalty. In a similar spirit, we formulate the general penalized estimation as
\[
\min_{\boldsymbol{f}} \frac{1}{2n} \sum_{i=1}^n \left( U_i(y) - \sum_{j=1}^p \langle f_j, \kappa_j(\cdot, X_{ij}) \rangle_{\mathcal{M}_j} \right)^2 + \sum_{j=1}^p P_{\lambda_1}\left( \|f_j\|_2 \right)+\frac{\lambda_2}{2} \sum_{j=1}^p \|f_j\|_{\mathcal{M}_j}^2,
\]

To simplify the notation, define $\mathcal{L}_n(\boldsymbol{f},\lambda_2) = \frac{1}{2n} \sum_{i=1}^n \left( U_i(y) - \sum_{j=1}^p \langle f_j, \kappa_j(\cdot, X_{ij}) \rangle_{\mathcal{M}_j} \right)^2+\frac{\lambda_2}{2} \sum_{j=1}^p \|f_j\|_{\mathcal{M}_j}^2$. We propose the following iterative algorithm.

\begin{algorithm}
\caption{The Local Linear Approximation (LLA) Algorithm for Fr\'echet Regression}\label{lla}
\begin{enumerate}
    \item \textbf{Initialization:} Start with  \( \hat{\boldsymbol{f}}^{(0)} = \boldsymbol{f}^{\text{init}} \) and 
    $
    \hat{w}_j^{(0)} = P'_{\lambda_1}(\|\hat{f}_j^{(0)}\|_2), \quad \text{for } j = 1, \dots, p.
    $
    
    \item \textbf{Iteration:} For \( m = 1, 2, \ldots \), repeat until convergence:
    \begin{enumerate}
        \item Update \( \hat{\boldsymbol{f}}^{(m)} \) by solving 
        $  
        \hat{\boldsymbol{f}}^{(m)} = \arg\min_{\boldsymbol{f}} \; \mathcal{L}_n(\boldsymbol{f},\lambda_2) + \sum_{j=1}^p \hat{w}_j^{(m-1)} \cdot \|f_j\|_2.
        $
        \item Update weights:
        $
        \hat{w}_j^{(m)} = P'_{\lambda_1}(\|\hat{f}_j^{(m)}\|_2), \quad \text{for } j = 1, \dots, p.
        $
    \end{enumerate}
\end{enumerate}
\end{algorithm}

Algorithm \ref{lla} requires an initial estimator \( \boldsymbol{f}^{(0)} \). In our implementation, we consider two options: (1) 
We apply the Elastic Net to a sequence of penalty values and select the optimal estimate minimizing the validation error, which is sparse and serves as a strong candidate for initial support recovery. (2) 
We employ ridge regression to obtain a smooth initial estimator. While ridge does not induce sparsity, the subsequent adaptive penalty can still achieve variable selection via appropriately chosen weights.

This framework effectively incorporates the strength of folded concave penalties while maintaining computational tractability. Section 4 provides theoretical guarantees showing that under suitable conditions, this procedure enjoys the \emph{strong oracle property}. This is motivated by the analogous results for SCAD in Euclidean regression \citep{fan2014strong}.

\section{Theoretical Properties}
\label{sec:theory}

This section first characterizes the optimality of solutions in Subsection 4.1, proves no false positives in variable selection in Subsection 4.2, and then establishes consistency and estimation accuracy in Subsection 4.3, as well as strong oracle property in Subsection 4.4.

\subsection{Optimality of Solutions}
\label{subsec:opt:sol}

We study the Karush–Kuhn–Tucker (KKT) conditions that characterize the optimality of the solution to the RKHS-based Elastic Net problem. 
Recall that
$
\mathcal{L}(f_1, \dots, f_p) := \frac{1}{2n} \sum_{i=1}^n \left( U_i(y) - \sum_{j=1}^p \langle f_j, \kappa_j(\cdot, X_{ij}) \rangle_{\mathcal{M}_j} \right)^2  + \lambda_1 \sum_{j=1}^p \|f_j\|_{\mathcal{M}_j} + \frac{\lambda_2}{2} \sum_{j=1}^p \|f_j\|_{\mathcal{M}_j}^2.
$ 

\begin{theorem}[KKT Conditions]\label{KKT-condition}
\(f_j \in \mathcal{M}_j\) is a solution to the optimization problem \eqref{optprob}, is equivalent to the following KKT conditions: 
\begin{itemize}
    \item[(i)] If \(\|f_j\|_{\mathcal{M}_j} > 0\), the coefficients \(\{\alpha_{ji}\}\) associated with \(f_j\) satisfy:   
    $$\frac{1}{n} \sum_{i=1}^n \left( U_i(y) - \sum_{k=1}^p \langle f_k, \kappa_k(\cdot, X_{ik}) \rangle_{\mathcal{M}_k} \right) \kappa_j(\cdot, X_{ij}) + \lambda_2 f_j = \lambda_1 \frac{f_j}{\|f_j\|_{\mathcal{M}_j}}.$$
    
    \item[(ii)] If \(\|f_j\|_{\mathcal{M}_j} = 0\), then:
    $$\left\| \frac{1}{n} \sum_{i=1}^n \left( U_i(y) - \sum_{k=1}^p \langle f_k, \kappa_k(\cdot, X_{ik}) \rangle_{\mathcal{M}_k} \right) \kappa_j(\cdot, X_{ij}) \right\|_{\mathcal{M}_j} \leq \lambda_1.$$
\end{itemize}
\end{theorem}

\begin{remark}
    If there exist two global minimizers, \(\boldsymbol{f}_1\) and \(\boldsymbol{f}_2\), for the optimization problem, then it must hold that: $\sum_{j=1}^p \|f_{1j}-f_{2j}\|_{\mathcal{M}_j}^2 = 0$. 
    This result follows directly from the KKT conditions' equivalence proof, which ensures the uniqueness of the global minimizer.
\end{remark}
\subsection{No False Positives in Variable Selection}
\label{subsec:noFP:var:sel}

We establish a partial consistency result for the proposed estimator. Specifically, we prove that with high probability, the set of selected variables is a subset of the true active set \( S_0 \). This ensures that no irrelevant variables are selected (i.e., no false positives), but the method does not guarantee inclusion of all truly relevant variables (i.e., false negatives may occur). 

Let the true active set be \( S_0 \subset \{1,\dots,p\} \) with \( s_0 = |S_0| \). We define the population covariance operator as $\Sigma = \mathbb{E} \left[ (\kappa_1(\cdot, X_1), \cdots, \kappa_p(\cdot, X_p)) \otimes (\kappa_1(\cdot, X_1), \cdots, \kappa_p(\cdot, X_p)) \right]$. One can estimate $\Sigma$ taking the expectation with respect to the empirical measure. 

We further define the
sub-operators \(\Sigma_{S_1, S_2}\) defined for any \(S_1, S_2 \subset \{1, \dots, p\}\). We define the \((\infty,\infty)\)-operator norm of operator \(\Sigma_{S_1 S_2}\) as
$
\|\Sigma_{S_1 S_2}\|_{\infty, \infty} = \sup_{\|f_1\|_\infty \leq 1,\, \|f_2\|_\infty \leq 1} \left\langle f_1,\, \Sigma_{S_1 S_2} f_2 \right\rangle,$
where \(f_1 = (f_j)_{j \in S_1}\) and \(f_2 = (f_k)_{k \in S_2}\) with each \(f_j \in \mathcal{M}_j\), and the supremum is taken over all such collections with $
\|f_1\|_\infty := \max_{j \in S_1} \|f_j\|_{\mathcal{M}_j}, \quad \|f_2\|_\infty := \max_{k \in S_2} \|f_k\|_{\mathcal{M}_k}.$
We use the same definition to denote the \((\infty,\infty)\)-norm for any other operator.

The following assumptions are made regarding the population-level covariance operator to ensure the correct identification of the true active set.
\begin{assumption}(Irrepresentable)
\label{ass:irrepresentable}
For some \(\gamma > 0\),   
   $ \| \Sigma_{S_0^cS_0} (\Sigma_{S_0S_0} + \lambda_2 I)^{-1} \|_{\infty,\infty} < 1 - \gamma.$    
\end{assumption}

\begin{assumption}(Scaling of Tuning Parameters)
\label{ass:scale:tuning}
For some constant \(C > 1\), the regularization parameters satisfy that $\frac{\lambda_1}{\lambda_2} > C \frac{1 - \gamma}{\gamma}$, $\frac{\sqrt{n} \lambda_2^2}{p} \to \infty$, and $\lambda_1 \to 0$.    
\end{assumption}

\begin{assumption}(Boundedness)
\label{ass:kernel}
There exists a uniform bound on the RKHS norm of the kernel feature map:
    $\sup_{j,x_j} \| \kappa_j(\cdot, x_j) \| < C$,
    and the signal norm, noise variance, and kernel expectations are uniformly bounded:
    $\|\boldsymbol{f}_0\| < C$, $\mathbb{E}[\epsilon^2] < C$, and $\sup_j \mathbb{E}[\kappa_j(X_j, X_j)] < C$.
    \end{assumption}

\begin{assumption}(Dimensionality)
\label{ass:dim}
The number of variables \(p = p_n\) may grow with \(n\), but must satisfy \( p = o(\sqrt{n}) \). The active set size \( s_0 \) is finite.
\end{assumption}
Assumption~\ref{ass:irrepresentable} is the irrepresentable condition, essential for variable selection consistency in Lasso, as it prevents excessive correlation between inactive and active variables \citep{zhao2006model, wainwright2009sharp}. Assumptions~\ref{ass:scale:tuning}--\ref{ass:dim} are standard in high-dimensional analysis \citep{van2008high, bickel2009simultaneous, negahban2012unified}. Specifically, Assumption~\ref{ass:scale:tuning} ensures proper scaling of regularization parameters to promote sparsity and stability; Assumption~\ref{ass:kernel} imposes uniform boundedness on kernel features, signal norms, and noise variance to enable concentration; and Assumption~\ref{ass:dim} controls dimensionality growth to preserve consistency. These conditions are widely adopted in high-dimensional regularized estimation and RKHS-based learning \citep{meier2009high, bickel2009simultaneous, koltchinskii2010oracle}.

\begin{theorem}[No False Positives]
\label{thm:subset}
Under the above assumptions, the estimated support set \(\hat{S}\) satisfies that 
$
\mathbb{P}(\hat{S} \subseteq S_0) \to 1.
$
\end{theorem}

The above theorem establishes a conservative property of the proposed estimator: irrelevant variables are excluded with high probability. 
The irrepresentable condition and tuning conditions used here are aligned with those in classical works such as \citep{zhao2006model, bickel2009simultaneous}, but the analysis is adapted to RKHS-based estimators.

\subsection{Variable Selection Consistency and Estimation Accuracy}
\label{subsec:accuracy}

We now establish that the proposed method achieves exact variable selection consistency, meaning that the estimated support \(\hat{S}\) recovers the true active set \(S_0\) with high probability. In addition, we provide uniform consistency of the estimated functions on the active coordinates.

\begin{assumption}(Signal Strength)
\label{ass:signal:str:active:set} 
    $\mathbb{E}[f_j^2(X_j)] > C$, for all $j \in S_0$, where \(C > 0\) is a constant. 
\end{assumption}

\begin{assumption}(Regularization Scaling)
\label{ass:regu:scaling}
\[
    \lambda_1\left\|\mathrm{diag}(\Sigma_{11},\Sigma_{22},\cdots,\Sigma_{s_0s_0})(\Sigma_{S_0S_0} + \lambda_2 I_{S_0})^{-1} \right\|_{\infty,\infty} = o(1).
    \]
\end{assumption}

Assumption~\ref{ass:signal:str:active:set} imposes a minimal signal strength condition on active components to avoid false negatives by ensuring that each active covariate contributes a sufficiently large signal relative to noise. Such assumptions are standard in high-dimensional variable selection to guarantee identifiability and consistency \citep{wainwright2009sharp, bickel2009simultaneous}. Assumption~\ref{ass:regu:scaling} controls the scaling of $\lambda_1$ relative to the covariance structure and dimensionality to prevent regularization from overwhelming the signal on the active set. While its form differs from prior work, it is not restrictive. For instance, \citet{zhao2006model} assumes a lower bound on the eigenvalues of $\Sigma_{11}$ and under $\lambda_2 \to 0$, our condition reduces to theirs. Similarly, \citet{zou2005regularization} requires the eigenvalues of the design covariance to be bounded away from zero and infinity, which implies our condition under suitable \(n,p\) scaling. More generally, such eigenvalue-based regularity is standard in high-dimensional regression; see \cite{wainwright2009sharp,bickel2009simultaneous,lee2016variable}, and recently \cite{liu2024robust}.

\begin{theorem}[Variable Selection and Estimation Consistency]
\label{thm:selection_consistency}
Under the above conditions and assumptions in Theorem \ref{thm:subset}, the proposed estimator satisfies:
\[
\mathbb{P}(\hat{S} = S_0) \to 1,
\quad \text{and} \quad
\sup_{j\in S_0}\|\Sigma_{jj}(\hat{f}_j-f_j) \| = o_p(1), \quad \text{for } j \in S_0.
\]
\end{theorem}

Theorem~\ref{thm:selection_consistency} establishes exact recovery of the sparsity pattern and accurate estimation of active effects. In comparison, \cite{tucker2023variable} only show that for non-active components \(j \notin S_0\), the estimated functions vanish in probability, i.e., \(\hat{f}_j \xrightarrow{p} 0\). This guarantee does not imply correct support recovery. Our result provides full control over both false positives and false negatives, leading to genuine variable selection consistency.

\subsection{Strong Oracle Property via LLA Algorithm}
\label{subsec:oracle:LLA}

We now analyze the sparsity recovery guarantees of the folded concave penalized Fr\'echet regression estimator obtained via the Local Linear Approximation (LLA) algorithm \citep{zou2008one,fan2014strong}. Let \( \boldsymbol{f} = (f_1, \dots, f_p) \) denote the additive component functions, and suppose the estimator solves the folded concave penalized problem:
\[
\min_{\boldsymbol{f}} \; \mathcal{L}_n(\boldsymbol{f},\lambda_2) + \sum_{j=1}^p P_\lambda(\|f_j\|_{\mathcal{M}_j}),
\]
where \( \mathcal{L}_n(\boldsymbol{f},\lambda_2) \) denotes the empirical Fr\'echet loss and \( \|\cdot\|_{\mathcal{M}_j} \) is the RKHS norm for the \(j\)-th coordinate. The penalty \( P_\lambda(\cdot) \) is assumed to satisfy the conditions (i)–(iv).

 Define the oracle estimator as
\[
\boldsymbol{f}^* = \arg\min_{\substack{f_j = 0,\ j \in S_0^c}} \mathcal{L}_n(\boldsymbol{f},\lambda_2).
\]
We analyze conditions under which the one-step LLA update recovers 
\( \boldsymbol{f}^* \) exactly. 
It is worth noting that our definition of the oracle estimator incorporates 
an additional $\ell_2$ penalty, since our problem is formulated in an infinite-dimensional RKHS. Without such $\ell_2$ penalty, the infinite-dimensionality of the RKHS would make the estimation ill-defined, even when the number of active predictors $p$ is small. In practice, however, this impact is negligible.

To this end, denote the subgradient of the empirical loss as \( \nabla \mathcal{L}_n(\boldsymbol{f},\lambda_2) \), and define the restricted gradient over \( S_0^c \) as $\nabla_{S_0^c} \mathcal{L}_n(\boldsymbol{f^*},\lambda_2) = \left\{ \nabla_j \mathcal{L}_n(\boldsymbol{f^*},\lambda_2) : j \in S_0^c \right\}$. Let $a_0=\min(1,a_2)$. 

\begin{theorem}[Strong Oracle Property via LLA]
\label{thm:lla-oracle}
Suppose the penalty function \( P_\lambda(\cdot) \) satisfies conditions (i)–(iv), and assume that 
$\min_{j \in S_0} \|f_j\|_{\mathcal{M}_j} > (a + 1)\lambda_1$. Let the initial estimator \( \boldsymbol{f}^{\mathrm{init}} \) be ridge-based or elastic-net-based, satisfying the conditions described below.

\begin{itemize}
    \item[\textbf{(a)}] \textbf{Ridge Initialization:} If
    $\|\boldsymbol{f}^{\mathrm{init}} - \boldsymbol{f}\|_{\infty} = o(\lambda_1)$, then, with high probability, the one-step LLA estimator satisfies that
    $\hat{\boldsymbol{f}}^{(1)} = \boldsymbol{f^*}$.

    \item[\textbf{(b)}] \textbf{Elastic Net Initialization:} If 
    $\|\boldsymbol{f}^{\mathrm{init}} - \boldsymbol{f}\|_{\max} \le a_0 \lambda_1$ and $\|\nabla_{A^c} \mathcal{L}_n(\boldsymbol{f}^{*})\|_{\max} < a_1 \lambda_1$ hold with high probability for \(a_0 + a_1 < a\) , then again,
    $\hat{\boldsymbol{f}}^{(1)} = \boldsymbol{f}^*$ with high probability.
\end{itemize}
\end{theorem}

Theorem~\ref{thm:lla-oracle} establishes a strong oracle property for the 
one-step LLA estimator under appropriate initialization. Our proof strategy 
builds on the framework of \cite{fan2014strong}. Related work includes \citet{lee2016variable}, who studied variable selection in RKHS via additive conditional independence but did not establish oracle properties, and \citet{zeng2014group}, who developed the SCAD-L$_2$ method in the classical finite-$p$ setting. Extending these results, we move beyond \cite{fan2014strong} to establish the strong oracle property for the SCAD-L$_2$ penalization in infinite-dimensional RKHS, thereby enabling applications to general metric-space-valued responses. To ensure Theorem~\ref{thm:lla-oracle} applies in practice, we also derive sufficient initialization accuracy for ridge- and elastic-net–based estimators.

\begin{lemma}[Ridge Initialization Accuracy]
\label{lem:ridge-init}
Assume the assumptions of the tuning parameter $\lambda_2$ and the rate of p in Theorem \ref{thm:selection_consistency} hold, then the ridge-based estimator satisfies that 
$\|\hat{\boldsymbol{f}}^{\mathrm{ridge}} - \boldsymbol{f}\|_{\max} = o(\lambda_2)$, with high probability as \( n \to \infty \).
\end{lemma}

\begin{lemma}[Elastic Net Initialization Accuracy]
\label{lem:enet-init}
Under assumptions of theorem \ref{thm:selection_consistency}, the Elastic Net estimator satisfies that $\|\hat{\boldsymbol{f}}^{\mathrm{enet}} - \boldsymbol{f}\|_{\max} \le a_0 \lambda_1$, and $\|\nabla_{S_0^c} \mathcal{L}_n(\boldsymbol{f^*},\lambda_2)\|_{\max} < a_1 \lambda_1$, with high probability as \( n \to \infty \).
\end{lemma}

\section{Estimation of Elastic Net in RKHS}
\label{sec:estimation}
Having established the theoretical foundation of the Elastic Net for variable selection in an additive Fr\'echet regression model, we now focus on the estimation strategy from the observed data $(Y_i, X_{ij}) \in (\Omega\lo Y, d\lo Y)\times \ca M\lo j,\ i=1,\dots,n;\ j= 1,\dots,p$. We first formally introduce Elastic Net in Reproducing Kernel Hilbert Spaces (RKHS) in a matrix representation. Using the coordinate representation and Gram matrix formalism, we express the functions \( f_j \) and the optimization problem in a compact matrix form.

The Representer Theorem is a foundational result in the theory of Reproducing Kernel Hilbert Spaces (RKHS), which ensures that the solution to a wide class of regularized empirical risk minimization problems admits a finite-dimensional representation. In the context of our additive Fr\'echet regression model, this theorem allows each functional component \( f_j \in \mathcal{M}_j \) to be expressed as a linear combination of kernel evaluations at the observed covariates, thus reducing the infinite-dimensional optimization to a finite-dimensional one, thus enabling an efficient computation by transforming the infinite-dimensional problem into one involving only \( n \) coefficients per component. 

Recall, from Section~\ref{subsec:additive:fr:model}, the evaluation of \( f_j \) at a sample point \( X_{ij} \) as $
f_j(X_{ij}) =$\\  $\sum_{k=1}^n \alpha_{jk} \kappa_j(X_{ij}, X_{kj}).$
Define the Gram matrix \( K_j \in \mathbb{R}^{n \times n} \) corresponding to the kernel \( \kappa_j \) as $
[K_j]_{ik} = \kappa_j(X_{ij}, X_{kj}), \quad i, k = 1, \dots, n.
$
Then, the evaluation of \( f_j \) at all sample points \( X_{ij} \) can be written in vectorized form as: $f_j(X_j) = K_j \alpha_j,$
where \( X_j = (X_{1j}, \dots, X_{nj})^\top \) is the vector of sample points in \( \Omega_{xj} \). Therefore, the RKHS norm of \( f_j \) is given by $
\|f_j\|_{\mathcal{M}_j}^2 = \alpha_j^\top K_j \alpha_j.$

Using the above representation, the Elastic Net problem in RKHS is now re-expressed  as 
\[
\min_{\alpha_1, \dots, \alpha_p} \frac{1}{2n} \|U(y) - \sum_{j=1}^p K_j \alpha_j \|_2^2 + \lambda_1 \sum_{j=1}^p \sqrt{\alpha_j^\top K_j \alpha_j} + \frac{\lambda_2}{2} \sum_{j=1}^p \alpha_j^\top K_j \alpha_j.
\]

Our proposed sparse Fr\'echet regression framework is implemented in reproducing kernel Hilbert spaces (RKHS) using two complementary optimization strategies that enable accurate and efficient variable selection: (i) an ADMM-based algorithm for solving the Elastic Net penalized least squares problem under RKHS norms, and (ii) an adaptive penalization scheme based on the right derivative of the SCAD penalty. The detailed derivation and algorithm for the ADMM updates are presented in Section~\ref{suppl:sec:admm} of the supplementary material. In practice, centering the Gram matrices and the response vector improves numerical stability and aligns with the assumption that the RKHS functions have zero mean.

\section{Simulation Studies}
\label{sec:simu}
In this section, we conduct extensive numerical studies to provide empirical evidence of the performance of our proposed variable selection framework in the context of global Fr\'echet regression. While Fr\'echet regression is inherently more abstract than standard linear regression due to the lack of vector space structure in the response domain, it remains highly applicable across a variety of modern data types. The implementation of the regression procedures depends on the choice of metric and the specific structure of the response space.

We consider two simulation scenarios for distributional data and symmetric positive definite (SPD) matrices, respectively. We analyze different data generation mechanisms for each type of random object --  one with additive linear effects and another with strong nonlinear signals in both mean and variance.  Across all simulations, the active set of covariates is predefined, allowing for rigorous assessment of variable selection accuracy.

We compare the following variable selection methods in our numerical experiments:
\begin{itemize}
    \item {Elastic Net Regularized Fr\'echet Regression (ElasticNet)}: our baseline method that applies Elastic Net penalization without the LLA step;
    \item {SCAD-L$_2$ Penalized Fr\'echet Regression with Ridge Initialization (RSCAD-L$_2$)}: which uses ridge regression for initialization, followed by SCAD-L$_2$ refinement;
    \item {SCAD-L$_2$ Penalized Fr\'echet Regression with Elastic Net Initialization (ESCAD-L$_2$)}: which uses the baseline Elastic Net as initialization for the SCAD-L$_2$ refinement.
\end{itemize}

All three proposed methods are compared against {Fr\'echet Ridge Selection Operator (FRiSO)} based on ridge-penalized global Fr\'echet regression \citep{tucker2023variable}.

Selection frequencies across repeated simulations are reported to assess both the \textit{accuracy} and \textit{sparsity} of the estimated models under varying signal complexities. Furthermore, to select the tuning parameter \( \lambda_1 \), we adopt a data-splitting strategy: the sample is divided into a training set and a test set. Model fitting is conducted on the training set, while prediction error is evaluated on the test set to determine the optimal tuning parameter.

Since kernel methods define function spaces via data-dependent Gram matrices, \( G_j^{\text{train}} \) and \( G_j^{\text{test}} \), the representations of estimated functions differ across datasets. To make the prediction consistent, we define the test-set evaluation as
\[
\hat{Y}_{\text{test},i} = \sum_{j=1}^p \sum_{k=1}^{n_{\text{train}}} \hat{\alpha}_{jk} \kappa_j(X^{\text{train}}_{kj}, X^{\text{test}}_{ij}),
\]
equivalently, using matrix notation,
$
\hat{\mathbf{Y}}_{\text{test}} = \sum_{j=1}^p \hat{\alpha}_j^\top \bar{G}_j,
$
where {\small \( (\bar{G}_j)_{k,i} =H_{train} \kappa_j(X^{\text{train}}_{kj}H_{test}, X^{\text{test}}_{ij}) \)}. The test-set mean squared error (MSE) is then defined as
\[
\text{MSE}_{\text{test}} = \| \widetilde{U}_{\text{test}} - \hat{\mathbf{Y}}_{\text{test}} \|^2,
\]
where \( \widetilde{U}_{\text{test}} = U(Y_{\text{test}}) - \bar{U}_{\text{train}} \) is centered using the training-set mean. This centering ensures compatibility between RKHS basis functions and the response.

\subsection{Distribution Data}
\label{subsec:simu:dist}
This section focuses on a regression model where the response \( Y \) is a univariate distribution characterized by its quantile function. We consider three data generation mechanisms, where the signal is incorporated following an additive model through location-scale transformation. 

\noindent\textbf{Model 1:}
Given covariates \( X = (X_1, \dots, X_p) \) with $p =30$, the conditional quantile function of \( Y \mid X \) is generated as:
\[
Y(t) = \mu_X + \sigma_X \Phi^{-1}(t), \quad t \in (0,1),
\]
with parameters: \( \mu_0 = 0 \), \( \beta = 3/4 \), \( \sigma_0 = 0 \), \( \gamma = 3 \), \( \nu_1 = 1 \), and \( \nu_2 = 0.5 \). Define
\[
\mu_X \sim \mathcal{N}(\mu_0 + \beta (X_4 + X_8), 1), \quad \sigma_X \sim \text{Gamma}\left(\frac{(\sigma_0 + X_1)^2}{v_2}, \frac{v_2}{\sigma_0 + X_1}\right).
\]
The covariates \( X \) are transformed from a multivariate Gaussian distribution with autoregressive covariance \(\Sigma_{ij} = \rho^{|i-j|}\), where \( \rho = 0.5 \), via \( X_j = 2\Phi(Z_j) - 1 \). By construction, the signal depends on \( X_1, X_4, X_8 \), which are thus treated as the active variables. We aim to recover this sparsity structure through variable selection performed using a linear kernel in the RKHS.

The results are summarized in the following table.
\begin{table}
\resizebox{.9\textwidth}{!}{
\begin{tabular}{lcccccccccc}
\toprule
Method & \textbf{X1} & X2 & X3 & \textbf{X4} & X5 & X6 & X7 & \textbf{X8} & X9 & X10 \\
\midrule
ElasticNet   & \textbf{0.99} & 0.07 & 0.13 & \textbf{1.00} & 0.04 & 0.02 & 0.10 & \textbf{1.00} & 0.01 & 0.01 \\
RSCAD-L$_2$ & \textbf{1.00} & 0.00 & 0.00 & \textbf{1.00} & 0.01 & 0.00 & 0.02 & \textbf{0.99} & 0.03 & 0.01 \\
ESCAD-L$_2$ & \textbf{1.00} & 0.05 & 0.08 & \textbf{1.00} & 0.02 & 0.04 & 0.06 & \textbf{1.00} & 0.04 & 0.06 \\
FRiSO      & \textbf{1.00} & 0.43 & 0.32 & \textbf{1.00} & 0.41 & 0.38 & 0.34 & \textbf{1.00} & 0.35 & 0.33 \\
\midrule\midrule
Method & X11 & X12 & X13 & X14 & X15 & X16 & X17 & X18 & X19 & X20 \\
\midrule
ElasticNet   & 0.01 & 0.02 & 0.01 & 0.02 & 0.01 & 0.01 & 0.01 & 0.01 & 0.02 & 0.01 \\
RSCAD-L$_2$ & 0.01 & 0.02 & 0.01 & 0.02 & 0.01 & 0.01 & 0.01 & 0.01 & 0.01 & 0.02 \\
ESCAD-L$_2$ & 0.04 & 0.03 & 0.07 & 0.02 & 0.04 & 0.05 & 0.04 & 0.08 & 0.05 & 0.07 \\
FRiSO      & 0.48 & 0.32 & 0.31 & 0.35 & 0.31 & 0.37 & 0.34 & 0.42 & 0.41 & 0.34 \\
\midrule\midrule
Method & X21 & X22 & X23 & X24 & X25 & X26 & X27 & X28 & X29 & X30 \\
\midrule
ElasticNet   & 0.01 & 0.03 & 0.00 & 0.02 & 0.00 & 0.00 & 0.01 & 0.03 & 0.02 & 0.01 \\
RSCAD-L$_2$ & 0.00 & 0.01 & 0.01 & 0.01 & 0.03 & 0.03 & 0.02 & 0.00 & 0.00 & 0.00 \\
ESCAD-L$_2$ & 0.01 & 0.08 & 0.02 & 0.03 & 0.02 & 0.02 & 0.05 & 0.09 & 0.05 & 0.05 \\
FRiSO      & 0.35 & 0.27 & 0.34 & 0.38 & 0.34 & 0.38 & 0.34 & 0.39 & 0.43 & 0.44 \\
\bottomrule
\end{tabular}
}
\centering
\caption{Variable selection frequencies under Model 1 with the linear distribution data and \( p = 30 \). The selection frequencies for active variables \textbf{X1},  \textbf{X4}  and \textbf{X8} are bolded.}
\end{table}

Compared to FRiSO, our proposed methods consistently demonstrate superior sparsity recovery. Notably, the {RSCAD-L$_2$} method achieves nearly perfect recovery of the true signals while maintaining low false positive rates. In contrast, FRiSO tends to over-select irrelevant variables despite successfully capturing the active ones. {ESCAD-L$_2$} provides a middle ground, offering good signal detection with slightly more noise inclusion. These results highlight the benefit of non-convex penalties and our RKHS-based framework in recovering the correct model structure under distribution-valued responses.

Next, we investigate a nonlinear setting where the true signal structure is no longer additive. To capture such nonlinearities, we employ a Gaussian kernel and examine the impact of kernel flexibility on variable selection accuracy.

We consider the following two nonlinear models to evaluate the performance of our proposed kernel-based methods. In both settings, the covariates \( X = (X_1, \dots, X_{10}) \) are generated from a multivariate normal distribution with autoregressive covariance \( \Sigma_{ij} = \rho^{|i-j|} \), transformed to the uniform scale via \( X_j = 2\Phi(Z_j) - 1 \), where \( \Phi \) is the standard Gaussian CDF. The response \( Y \) is the distributional object:$
Y(\cdot) = \mu_X + \sigma_X \Phi^{-1}(\cdot),
$ this time with a nonlinear signal embedded in \( \mu_X \) and \( \sigma_X \).

\noindent \textbf{Model 2.}  
This model introduces an exponential nonlinear structure:
\[
\mu_X \sim \mathcal{N} \left( \mu_0 + \beta \left( \exp(-X_4^2) + \exp(-X_8^2) \right), \nu_1 \right),
\]
\[
\sigma_X \sim \text{Gamma} \left( \frac{(\sigma_0 + \gamma \exp(-2(X_1 - 1)^2))^2}{\nu_2}, \frac{\nu_2}{\sigma_0 + \gamma \exp(-2(X_1 - 1)^2)} \right),
\]
with parameters: \( \mu_0 = 0 \), \( \beta = 12 \), \( \sigma_0 = 0 \), \( \gamma = 12 \), \( \nu_1 = 1 \), and \( \nu_2 = 0.5 \).

\noindent \textbf{Model 3.}  
This model contains trigonometric and rational polynomial nonlinearities:
\[
\mu_X \sim \mathcal{N} \left( \mu_0 + \beta \left( \sin(2\pi X_4) + \frac{2}{1 + |X_8|} \right), v_1 \right),
\]
\[
\sigma_X \sim \text{Gamma} \left( \frac{(\sigma_0 + \gamma \exp(-(X_1 - 1)^2))^2}{v_2}, \frac{v_2}{\sigma_0 + \gamma \exp(-(X_1 - 1)^2)} \right),
\]
with parameters: \( \mu_0 = 0 \), \( \beta = 10 \), \( \sigma_0 = 0 \), \( \gamma = 20 \), \( v_1 = 1 \), and \( v_2 = 0.5 \).

In both models, the true signal variables are \( \textbf{X1}, \textbf{X4}, \textbf{X8} \), and only our Gaussian-kernel-based methods accurately recover all three signal variables with almost no false positives. In contrast, the linear-kernel versions struggle due to their limited capacity to capture nonlinear dependencies. FRiSO, while detecting the signals, suffers from substantial over-selection.
The results are summarized in the following table, and the solution path corresponding to the elastic net results is provided in the Section~\ref{suppl:sec:simu:sol:path}

 \begin{table}
\resizebox{.9\textwidth}{!}{
\begin{tabular}{lcccccccccc}
\toprule
Method & \textbf{X1} & X2 & X3 & \textbf{X4} & X5 & X6 & X7 & \textbf{X8} & X9 & X10 \\
\midrule
ElasticNet (Linear)       & \textbf{0.47} & 0.09 & 0.10 & \textbf{0.08} & 0.04 & 0.10 & 0.11 & \textbf{0.06} & 0.07 & 0.08 \\
ElasticNet (Gaussian)     & \textbf{1.00} & 0.00 & 0.00 & \textbf{1.00} & 0.01 & 0.00 & 0.00 & \textbf{1.00} & 0.01 & 0.00 \\
RSCAD-L$_2$ (Linear)     & \textbf{0.39} & 0.13 & 0.09 & \textbf{0.12} & 0.11 & 0.12 & 0.10 & \textbf{0.09} & 0.10 & 0.12 \\
RSCAD-L$_2$ (Gaussian)   & \textbf{1.00} & 0.14 & 0.13 & \textbf{1.00} & 0.14 & 0.20 & 0.14 & \textbf{1.00} & 0.16 & 0.16 \\
ESCAD-L$_2$ (Linear)     & \textbf{0.20} & 0.06 & 0.05 & \textbf{0.03} & 0.02 & 0.04 & 0.07 & \textbf{0.02} & 0.03 & 0.02 \\
ESCAD-L$_2$ (Gaussian)   & \textbf{1.00} & 0.01 & 0.00 & \textbf{1.00} & 0.00 & 0.00 & 0.00 & \textbf{1.00} & 0.00 & 0.00 \\
FRiSO                   & \textbf{1.00} & 0.52 & 0.51 & \textbf{0.41} & 0.49 & 0.49 & 0.57 & \textbf{0.52} & 0.50 & 0.53 \\
FRiSO (Refit)           & \textbf{1.00} & 0.09 & 0.10 & \textbf{0.08} & 0.06 & 0.09 & 0.06 & \textbf{0.09} & 0.06 & 0.07 \\
\bottomrule
\end{tabular}}
\centering
\caption{Variable selection frequencies under Model 2 with the nonlinear distribution data and \( p = 10 \). The selection frequencies for active variables \textbf{X1},  \textbf{X4} and \textbf{X8} are bolded.}
\label{tab:nonlinear_model1}
\end{table}

\begin{table}
\resizebox{.9\textwidth}{!}{\begin{tabular}{lcccccccccc}
\toprule
Method & \textbf{X1} & X2 & X3 & \textbf{X4} & X5 & X6 & X7 & \textbf{X8} & X9 & X10 \\
\midrule
ElasticNet (Linear)       & \textbf{0.98} & 0.03 & 0.04 & \textbf{0.02} & 0.05 & 0.03 & 0.03 & \textbf{0.04} & 0.07 & 0.03 \\
ElasticNet (Gaussian)     & \textbf{1.00} & 0.04 & 0.06 & \textbf{1.00} & 0.15 & 0.03 & 0.04 & \textbf{1.00} & 0.02 & 0.04 \\
RSCAD-L$_2$ (Linear)     & \textbf{0.95} & 0.11 & 0.03 & \textbf{1.00} & 0.09 & 0.05 & 0.08 & \textbf{0.08} & 0.06 & 0.07 \\
RSCAD-L$_2$ (Gaussian)   & \textbf{1.00} & 0.16 & 0.21 & \textbf{1.00} & 0.19 & 0.10 & 0.17 & \textbf{1.00} & 0.13 & 0.11 \\
ESCAD-L$_2$ (Linear)     & \textbf{0.96} & 0.02 & 0.00 & \textbf{0.01} & 0.01 & 0.01 & 0.00 & \textbf{0.00} & 0.01 & 0.01 \\
ESCAD-L$_2$ (Gaussian)   & \textbf{0.99} & 0.05 & 0.01 & \textbf{1.00} & 0.01 & 0.00 & 0.01 & \textbf{1.00} & 0.00 & 0.01 \\
FRiSO                   & \textbf{1.00} & 0.28 & 0.50 & \textbf{1.00} & 0.44 & 0.40 & 0.45 & \textbf{0.40} & 0.43 & 0.44 \\
FRiSO (Refit)           & \textbf{1.00} & 0.03 & 0.15 & \textbf{1.00} & 0.10 & 0.07 & 0.08 & \textbf{0.06} & 0.03 & 0.06 \\
\bottomrule
\end{tabular}
}
\centering
\caption{Variable selection frequencies under Model 3 with the nonlinear distribution data and \( p = 10 \). The selection frequencies for active variables \textbf{X1} \textbf{X4}, and \textbf{X8} are bolded.}
\label{tab:nonlinear_model2}
\end{table}

Similar to Model 1, Gaussian-kernel methods show significantly improved performance over their linear counterparts. In particular, RSCAD-L$_2$ and ESCAD-L$_2$ under Gaussian kernels reliably select all signal variables with very limited noise inclusion. In contrast, FRiSO tends to over-select and lacks stability under strong nonlinearity.

These experiments confirm the critical role of flexible kernels for recovering nonlinear relationships and demonstrate the practical advantage of our framework in such scenarios.

\subsection{SPD Matrix Data}
\label{subsec:simu:spd}

We next consider object responses that lie in the space of symmetric positive definite (SPD) matrices. In this setting, the signal is encoded through an additive structure in the mean of the Cholesky factor of the SPD matrix via both a linear and nonlinear structural dependence. Variable selection is performed under a non-Euclidean geometry induced by the Cholesky decomposition metric.  Let \( X = (X_1, \ldots, X_p) \) be a \( p \)-dimensional vector of covariates drawn from a multivariate normal distribution with an autoregressive correlation structure:
\[
X \sim \mathcal{N}(0, \Sigma), \quad \Sigma_{ij} = \rho^{|i - j|}, \quad \rho = 0.5.
\]

\noindent \textbf{Model 4:}
We generate SPD responses through a latent Cholesky factor \( A \), whose conditional expectation depends additively on a sparse subset of covariates:
\[
E(A \mid X) = (\mu_0 + \beta (X_1 + X_3)) I + (\sigma_0 + \gamma (X_5 + X_7 + X_9)) U.
\]
Here, \( I \) is the identity matrix, and \( U \) is an upper triangular matrix with ones above the diagonal. The observed response is constructed as
$
Y = A^\top A,
$ which lies in the space of SPD matrices denoted \( \mathcal{TP}_1^2 \). The metric used to evaluate distances between SPD matrices is the squared Frobenius norm of the Cholesky factors:
\[
d_C(P_1, P_2) = \| P_1^{1/2} - P_2^{1/2} \|_F^2.
\]
This setup induces a sparse regression model with active variables \( X_1, X_3, X_5, X_7, X_9 \). The goal is to identify this active set using various variable selection procedures; we employ linear kernels to capture the additive structure of the covariates. Table~\ref{tab:spd_results_clean} reports variable selection frequencies across 100 simulation replicates under a dimensionality of \( p = 30 \). Again, we compare the proposed methods—ElasticNet, RSCAD-L$_2$, and ESCAD-L$_2$ with the FRiSO.
\begin{table}
\resizebox{.9\textwidth}{!}{
\resizebox{.9\textwidth}{!}{
\begin{tabular}{lcccccccccc}
\toprule
Method & \textbf{X1} & X2 & \textbf{X3} & X4 & \textbf{X5} & X6 & \textbf{X7} & X8 & \textbf{X9} & X10 \\
\midrule
ElasticNet     & \textbf{0.96} & 0.14 & \textbf{0.97} & 0.25 & \textbf{1.00} & 0.14 & \textbf{1.00} & 0.07 & \textbf{1.00} & 0.07 \\
RSCAD-L$_2$   & \textbf{1.00} & 0.00 & \textbf{0.96} & 0.02 & \textbf{1.00} & 0.02 & \textbf{1.00} & 0.00 & \textbf{1.00} & 0.00 \\
ESCAD-L$_2$   & \textbf{0.81} & 0.02 & \textbf{0.84} & 0.03 & \textbf{1.00} & 0.01 & \textbf{1.00} & 0.01 & \textbf{1.00} & 0.02 \\
FRiSO        & \textbf{1.00} & 0.00 & \textbf{1.00} & 0.02 & \textbf{1.00} & 0.01 & \textbf{1.00} & 0.00 & \textbf{1.00} & 0.00 \\
\midrule\midrule
Method & X11 & X12 & X13 & X14 & X15 & X16 & X17 & X18 & X19 & X20 \\
\midrule
ElasticNet     & 0.00 & 0.04 & 0.01 & 0.02 & 0.03 & 0.03 & 0.01 & 0.01 & 0.01 & 0.01 \\
RSCAD-L$_2$   & 0.02 & 0.02 & 0.02 & 0.04 & 0.02 & 0.02 & 0.02 & 0.04 & 0.06 & 0.02 \\
ESCAD-L$_2$   & 0.00 & 0.00 & 0.01 & 0.02 & 0.02 & 0.01 & 0.02 & 0.02 & 0.02 & 0.00 \\
FRiSO        & 0.00 & 0.00 & 0.00 & 0.00 & 0.00 & 0.00 & 0.00 & 0.00 & 0.00 & 0.00 \\
\midrule\midrule
Method & X21 & X22 & X23 & X24 & X25 & X26 & X27 & X28 & X29 & X30 \\
\midrule
ElasticNet     & 0.02 & 0.02 & 0.01 & 0.05 & 0.00 & 0.03 & 0.01 & 0.00 & 0.02 & 0.00 \\
RSCAD-L$_2$   & 0.00 & 0.08 & 0.02 & 0.02 & 0.02 & 0.00 & 0.02 & 0.02 & 0.00 & 0.04 \\
ESCAD-L$_2$   & 0.01 & 0.01 & 0.01 & 0.02 & 0.01 & 0.00 & 0.00 & 0.00 & 0.00 & 0.00 \\
FRiSO        & 0.00 & 0.00 & 0.00 & 0.00 & 0.00 & 0.00 & 0.00 & 0.00 & 0.00 & 0.00 \\
\bottomrule
\end{tabular}
}
}
\centering
\caption{Variable selection frequencies under Model 4 with the linear SPD matrix and \( p = 30 \). The selection frequencies for active variables  \textbf{X1}, \textbf{X3}, \textbf{X5}, \textbf{X7}, and \textbf{X9} are bolded.}
\label{tab:spd_results_clean}
\end{table}

The {RSCAD-L$_2$} procedure exhibits the strongest performance, consistently identifying all true active variables while maintaining near-zero false discovery. {ESCAD-L$_2$} also recovers the signal effectively, with slightly more inclusion of irrelevant features. The {ElasticNet} method performs reasonably well but tends to over-select mid-indexed covariates, reducing overall sparsity. The {FRiSO} method correctly identifies the active set, yet its stability and specificity are sensitive to the post-thresholding mechanism, as further illustrated in Section~\ref{suppl:sec:friso_refit} of the supplements. These results highlight the effectiveness of our RKHS-based procedures with non-convex penalties for model selection in manifold-valued responses, particularly under structured linear signal models on SPD spaces.

We further evaluate our methods under a nonlinear SPD regression model to assess robustness against nonlinear mean and variance structures. In this setting, the response remains an SPD matrix observed as \( Y = A^\top A \), while both the mean and scale parameters of the underlying matrix \( A \) depend nonlinearly on the covariates.\\

\noindent \textbf{Model 5:} The nonlinear structure is introduced as follows
\begin{align*}
\mu_X &= \mu_0 + \beta \left(3X_1^2 + \sin(2\pi X_3)\right), \text{ and }
\sigma_X \sim \text{Gamma} \left( \frac{(\sigma_0 + \gamma \eta_X)^2}{v_2}, \, \frac{v_2}{\sigma_0 + \gamma \eta_X} \right),
\end{align*}
where $\eta_X = \exp(-X_5) + 2\exp(-2(X_7 - 1)^2) + \frac{2}{1 + |X_9|},$ with \( \mu_0 = 3,  \beta = 2, \sigma_0 = 1 \), \( \gamma = 4 \), \( v_1 = 1 \), and \( v_2 = 0.5 \). The Cholesky distance is used to compare SPD matrices.

The active covariates under this setting are \( \textbf{X1}, \textbf{X3}, \textbf{X5}, \textbf{X7}, \textbf{X9} \). Table~\ref{tab:spd_nonlinear_results} reports the selection frequencies under the nonlinear SPD model with \( n = 200 \) and \( p = 10 \). Our proposed methods—particularly RSCAD-L$_2$ and ESCAD-L$_2$—demonstrate strong performance in recovering the correct support: all five active variables (\textbf{X1}, \textbf{X3}, \textbf{X5}, \textbf{X7}, \textbf{X9}) are selected with high frequency, and false positives are largely suppressed.

In contrast, FRiSO exhibits considerable instability in this nonlinear regime. While it successfully includes some signal variables such as \textbf{X1}, \textbf{X5}, and \textbf{X9}, it frequently misses \textbf{X3} and \textbf{X7}, and selects several noise variables. This instability is further exacerbated in the post-thresholded version, FRiSO (Refit), which becomes overly conservative and fails to retain key active variables. These observations reaffirm the sensitivity of FRiSO to nonlinear structure and thresholding choices and highlight the robustness of our kernel-based methods with non-convex regularization in challenging settings.

\begin{table}
\begin{tabular}{lcccccccccc}
\toprule
Method & \textbf{X1} & X2 & \textbf{X3} & X4 & \textbf{X5} & X6 & \textbf{X7} & X8 & \textbf{X9} & X10 \\
\midrule
ElasticNet       & \textbf{1.00} & 0.01 & \textbf{0.74} & 0.30 & \textbf{1.00} & 0.06 & \textbf{1.00} & 0.06 & \textbf{1.00} & 0.19 \\
RSCAD-L$_2$     & \textbf{1.00} & 0.16 & \textbf{1.00} & 0.08 & \textbf{1.00} & 0.06 & \textbf{1.00} & 0.10 & \textbf{1.00} & 0.08 \\
ESCAD-L$_2$     & \textbf{1.00} & 0.01 & \textbf{0.82} & 0.04 & \textbf{0.99} & 0.02 & \textbf{1.00} & 0.00 & \textbf{1.00} & 0.01 \\
FRiSO          & \textbf{1.00} & 0.05 & \textbf{0.02}         & 0.11 & \textbf{1.00} & 0.05 & \textbf{0.03}         & 0.01 & \textbf{1.00} & 0.00 \\
FRiSO (Refit)  & \textbf{1.00} & 0.00 & \textbf{0.01}         & 0.00 & \textbf{1.00} & 0.03 & \textbf{0.02}         & 0.00 & \textbf{1.00} & 0.00 \\
\bottomrule
\end{tabular}
\centering
\caption{Variable selection frequencies under Model 5 with the nonlinear SPD matrix and \( p = 10 \). The selection frequencies for active variables \textbf{X1}, \textbf{X3}, \textbf{X5}, \textbf{X7}, and \textbf{X9} are bolded.}
\label{tab:spd_nonlinear_results}
\end{table}

\subsection{Empirical Validation of Sparsity Invariance}
\label{subsec:simu:sparsity:inv}
The stability of the active set, offered by the sparsity invariance assumption, which is fundamental to our variable selection method, is further tested through numerical studies. We make a careful investigation on the robustness of variable selection results in the context of both distribution-valued and SPD matrix-valued responses in the Fr\'echet regression framework, under several perturbed choices of \( y \), and report the selection frequencies of individual covariates across different models and sample sizes. The results empirically confirm that, while the signal strength may vary with \( y \), the set of selected variables remains stable and consistently converges to the true support as the sample size grows. More details are presented in Section~\ref{suppl:subsec:simu:sparsity:inv} of the Supplements.

\section{Application to Bike Rental Data}
\label{sec:data}

\begin{figure}[ht]
\centering
\begin{minipage}{.5\textwidth}
  \centering
  \includegraphics[width=.9\textwidth]{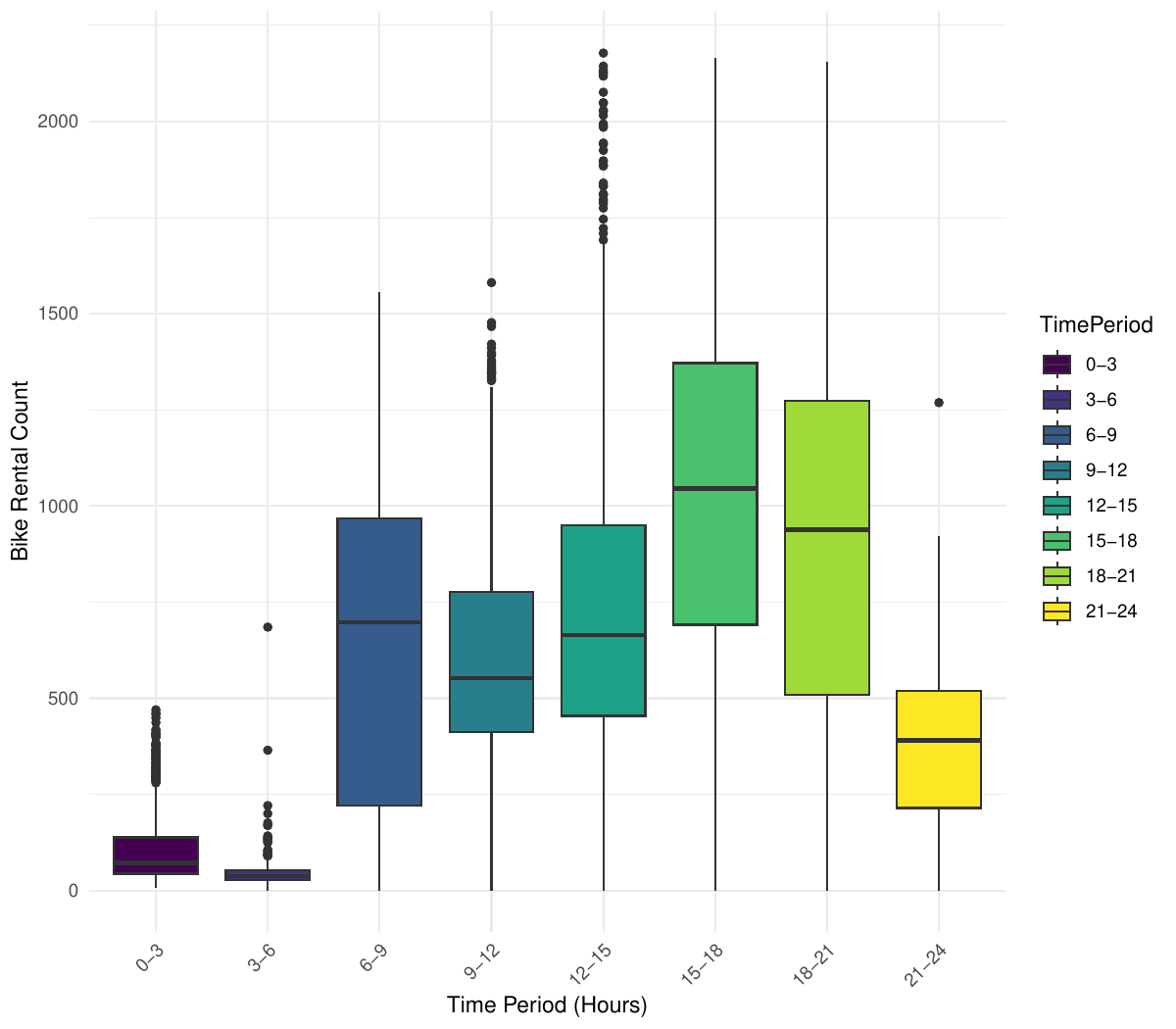}
\end{minipage}%
\begin{minipage}{.5\textwidth}
  \centering
  \includegraphics[width=.9\textwidth]{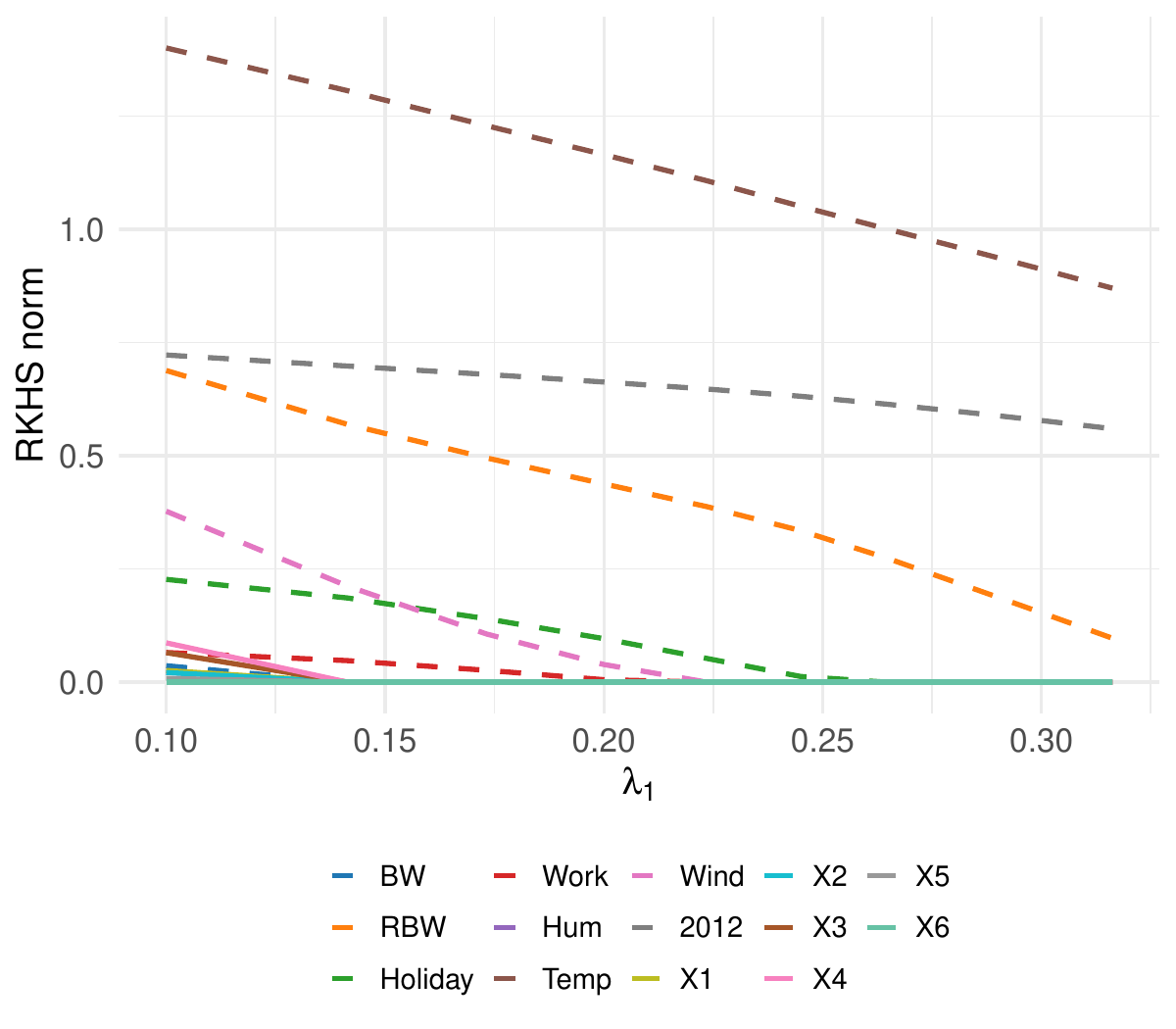}
\end{minipage}
\centering
\caption{The left panel shows the boxplot of hourly bike rentals over all days. In the right panel, the solution paths of all covariates under the Elastic Net method, plotted across a grid of $\lambda_1$ values, showcase the convergence of the proposed optimization algorithm.}
\label{fig:bike:rental}
\end{figure}

We analyze real-world data from the Capital Bikeshare system in Washington, D.C., consisting of hourly rental records from 2011–2012. Each day provides 24 hourly observations, yielding a total of 731 daily rental profiles. The response \( Y_i \in \mathbb{R}^{24} \) denotes the cumulative bike rentals up to each hour for day \( i \). From this, we construct the response for each day as a 24-dimensional vector representing the empirical quantiles of bike rental counts. Thus, each observation is treated as a \textit{distribution-valued object} residing in the space of univariate distributions, equipped with the Wasserstein-2 metric, which quantifies the distance between two distributional objects with cumulative distribution functions $H(\cdot)$ and $G(\cdot)$ as:
\[
W_2^2(H, G) = \int_0^1 \left[ H^{-1}(t) - G^{-1}(t) \right]^2 dt.
\]
The left panel of Figure~\ref{fig:bike:rental} provides a visualization of the rental pattern across 731 days. Clear peaks during morning and evening commuting hours motivate a functional approach.

The covariates \( X_i \in \mathbb{R}^p \) include eight meaningful environmental and calendar-based covariates, as well as synthetic noise variables.
The core covariates include: BW, RBW: Indicators of bad or severe weather (e.g., mist, rain, snow); Holiday, Working Day: Calendar-based indicators;  Humidity, Temperature, Windspeed: Daily weather statistics;  2012: A year indicator.
All covariates are standardized. 

We evaluate our three proposed procedures (ElasticNet, RSCAD-L$_2$, ESCAD-L$_2$) and compare them with the variable selection method, FRiSO, as proposed in~\cite{tucker2023variable}. FRiSO applies a ridge-penalized global Fr\'echet regression and selects variables via a post-thresholding step. In bike rental analysis, FRiSO performs variable selection based on the fitted solution path, with an additional refitting step to enhance interpretability by excluding weak covariates. Two scenarios for variable selection are considered to compare the performance of these methods: (a) {Moderate-dimensional setting} with 8 observed covariates + 6 synthetic noise variables (\(p=14\)); (b) {High-dimensional setting} with 8 observed covariates + 24 synthetic noise variables (\(p=32\)).

\vspace{0.5em}
\noindent\textbf{(a) Moderate-dimensional setting (\(p=14\)).}

Table~\ref{tab:bike_p14_reformatted} shows the variable selection results for setting(a). All three of our methods successfully eliminated the noise variables. In contrast, FRiSO selected multiple synthetic covariates (X2, X5), raising concerns about overfitting. 
To further illustrate the variable selection behavior, we take the standard method as an example and visualize the solution paths of all covariates across a range of $\lambda_1$ values. As shown in the right panel of Figure~\ref{fig:bike:rental}, the results are consistent with the selection results reported in Table~\ref{tab:bike_p14_reformatted}. All noise variables, as well as two covariates---BW and Hum---which are considered unimportant by our proposed method, exhibit rapid shrinkage and converge to zero as $\lambda_1$ increases.

\begin{table}
\begin{tabular}{l|c}
\toprule
Method & Selected Covariates \\
\midrule
FRiSO       & BW, RBW, Holiday, Work, Hum, Temp, Wind, 2012, X2, X5\\
ElasticNet    & RBW, Holiday, Work, Temp, Wind, 2012\\
RSCAD-L$_2$  & RBW, Holiday, Work, Temp, Wind, 2012\\
ESCAD-L$_2$  & RBW, Holiday, Work, Temp, Wind, 2012\\
\bottomrule
\end{tabular}
\centering
\caption{Variable selection results with 10-fold CV under the moderate-dimensional setting.}
\label{tab:bike_p14_reformatted}
\end{table}

\vspace{1em}
\noindent\textbf{(b) High-dimensional setting (\(p=32\)).}
\begin{table}
\begin{tabular}{l|c}
\toprule
{Method} & {Selected Covariates} \\
\midrule
FRiSO       & BW, RBW, Holiday, Work, Hum, Temp, Wind, 2012, X2, X6, X8, X14 \\
ElasticNet    & RBW, Holiday, Work, Temp, Wind, 2012, X3, X6\\
RSCAD-L$_2$  & RBW, Holiday, Work, Temp, Wind, 2012\\
ESCAD-L$_2$  & RBW, Holiday, Work, Temp, Wind, 2012\\
\bottomrule
\end{tabular}
\centering
\caption{Variable selection results with 10-fold CV under the high-dimensional setting. 
}
\label{tab:bike_p32_final}
\end{table}

Table~\ref{tab:bike_p32_final} showcases the superiority of our approaches for variable selection in the higher-dimensional regime. Our methods exhibit strong performance in both moderate- and high-dimensional regimes. In the moderate-dimensional setting (\(p=14\)), all three proposed procedures successfully eliminated the six synthetic noise variables while retaining meaningful covariates such as Holiday, Work, and RBW. In contrast, FRiSO selected two synthetic variables (X2, X5), indicating a tendency toward overfitting due to its lack of sparsity control.

In the high-dimensional setting (\(p=32\)), the advantage of our methods becomes clearer. FRiSO selected seven irrelevant synthetic variables across different blocks, whereas RSCAD-L$_2$ and ESCAD-L$_2$ achieved perfect exclusion of all 24 noise variables. This highlights the robustness of our penalized framework under high-dimensional noise contamination.

Moreover, our methods consistently selected covariates that are both interpretable and directly relevant to bike rental behavior. For example, Holiday and Work reflect calendar-based effects that naturally influence rental data, while RBW (really bad weather) captures extreme conditions such as rain or snow that are known to deter bike usage. In contrast, variables like BW (mildly bad weather) and Humidity were not selected by our methods. This suggests that their effects are either negligible or largely captured by stronger, more informative covariates such as RBW and Temperature. Also, FRiSO applies a refitting step to improve interpretability, but this adjustment is not guided by a sparsity-inducing criterion and lack a theoretical guarantee. It may compromise stability in high dimensions.

To conclude, these results demonstrate that our SCAD-based methods provide more stable, sparse, and interpretable variable selection results, especially in the presence of noise and increased dimensionality.

\section{Conclusion}
\label{sec:concl}
In this paper, we introduced a novel additive Fréchet regression model for variable selection in object-valued data analysis, addressing a key gap in the literature. Existing Fréchet regression methods lack principled mechanisms for sparsity. By transforming squared-distance responses into additive components and applying an elastic-net-type penalty within an RKHS framework, we developed a regularized estimator that is both geometrically valid and computationally efficient. Our theoretical results establish selection consistency, and empirical evidence shows improved predictive performance over existing non-Euclidean models.

After identifying relevant covariates, a natural next step is to refit the model using only selected variables. Here, nonlinear global Fréchet regression methods—such as the RKHS-based estimators of \citet{bhat:25}—are especially useful, offering optimal convergence rates under mild assumptions. Further investigation into the asymptotics of the refitted model could clarify trade-offs between parsimony and accuracy. We conjecture that incorporating debiasing techniques, similar to those used in desparsified Elastic Net approaches in Euclidean settings, may yield sharper convergence rates and facilitate valid post-selection inference, including confidence sets for effect sizes.

To summarize, our primary contribution lies in bridging an important methodological gap by providing a theoretically grounded and scalable approach to variable selection in Fréchet regression. Unlike previous methods, which focused solely on estimation or prediction, our framework integrates nonlinear additive modeling with elastic-net-style regularization to enable both sparse recovery and smooth estimation. Furthermore, our idea of targeting the squared metric $d \hi 2(Y, y)$ rather than the random object  $Y$ itself when performing variable selection allows us to avoid imposing any additional structure on $Y$ beyond that it is a member of a metric space.  This means our method is applicable to {\em all} metric-space-valued responses.
To our knowledge, this is the first work to establish variable selection theory for regression with general metric-space-valued responses. Simulations confirm that our method recovers true active covariates and improves prediction, underscoring its practical value and broad applicability to complex, non-Euclidean data.

{\small
\renewcommand{\baselinestretch}{1}

\bibliographystyle{agsm}
\bibliography{ref}
}
\end{document}